\def\d{\bar{\partial}}

\magnification\magstep1   
\baselineskip=18pt
\let\A=\i

\def\i{\infty}
\def\l{\ell}

\def\ra{\rightarrow}
\def\a{\langle} \def\b{\rangle}
\def\w{\widetilde}
\def\i{\infty}
\def\v{\vert}
\def\V{\Vert}
\def\n{\noindent}
\centerline{{\bf  Interpolation Between
$H^p$ Spaces and Non-Commutative
Generalizations II}\footnote*{Supported in part by N.S.F.
grant DMS 9003550}} \vskip12pt \centerline {by  Gilles
Pisier} \vskip12pt

\centerline{\bf Abstract}  
We continue an investigation started in a preceding
paper. We discuss the classical results of Carleson
connecting Carleson measures with the  $\d$-equation in a
slightly more abstract framework than usual. We also
consider  a more recent result of Peter Jones which
shows the  existence of a solution of the $\d$-equation,
which satisfies simultaneously a good $L_\i$ estimate and
a good $L_1$ estimate. This appears as a special case of
our main result which can be stated as follows: Let
 $(\Omega,\cal{A},\mu)$ be  
any measure space. Consider a bounded operator $u:H^1\ra
L_1(\mu)$. Assume that on one hand $u$ admits an extension
$u_1:L^1\ra L_1(\mu)$  bounded with norm $C_1$, and on
the other hand that $u$ admits an extension
$u_\i:L^\i\ra L_\i(\mu)$  bounded with norm $C_\i$. Then 
$u$ admits  an extension $\w{u}$ which is bounded
simultaneously from $L^1$ into $ L_1(\mu)$ and from $L^\i$
into $ L_\i(\mu)$ and satisfies $$\eqalign{&\|\tilde
u\colon \ L_\infty \to L_\infty(\mu)\|\le CC_\infty\cr
&\|\tilde u\colon \ L_1\to L_1(\mu)\|\le CC_1}$$ where $C$
is a numerical constant. \vskip6pt
\n{\bf Introduction}  
\medskip

We will denote by  $D$ the open unit
disc of the complex plane, by $\bf{T}$   the unit
circle and by $m$ the normalized Lebesgue
measure on $\bf{T}$.
Let $0<p\leq\i$. We will denote simply by $L_p$ the
space   $L_p({\bf T},m)$ and
by $H^p$ the classical Hardy space of
analytic functions on   $D$. It is well known that
$H^p$ can be identified with a closed subspace of  $L_p$,
namely the closure in $L_p$ (for $p=\i$ we must take the
weak*-closure) of the linear span of the functions
$\{e^{int}|n\ge 0\}$. More generally, when $B$ is a Banach
space, we denote by $L_p(B)$ the  usual space of
Bochner-$p$-integrable $B$-valued functions on $({\bf 
T},m)$, so that when $p<\i$, $L_p \otimes B$ is dense in
$L_p(B)$. We denote by  
$H^p(B)$ (and simply $H^p$ if $B$ is one dimensional)  
the Hardy space of B-valued analytic functions $f$ such
that  $\sup_{r<1}(\int\|f(rz)\|^p dm(z))^{1/p} <\i.$
 We denote $$\V f\V_{H^p(B)} =
\sup_{r<1}(\int\|f(rz)\|^p dm(z))^{1/p}.$$
We refer   to [G] and [GR]
for more information on $H^p$-spaces and to [BS] and [BL]
for more on real and complex interpolation.

We recall that a finite positive
measure $\mu$ on $D$   is called a Carleson measure if
there is a constant $C$ such that for any $r>0$ and any
real number $\theta$, we have $$\mu(\{z\in D, \ \v
z\v>1-r ,\ \v \arg(z)- {\theta}\v <r\}) \le Cr.$$ We will
denote by $\V \mu \V_C$ the smallest constant $C$ for
which this holds. Carleson (see [G]) proved that, for
each $0<p<\i$, this norm  $\V \mu \V_C$ is equivalent to
the smallest constant $C'$ such that  $$\forall f\in
H^p\qquad \int \v f\v^p d\mu \le C' \V f\V_{H^p}^p
.\leqno(0.1)$$ Moreover, he proved that, for any $p>1$ there
is a constant $A_p$ such that any harmonic function $v$ on
$D$ admitting radial limits in $L_p({\bf T},m)$  
satisfies
$$  \int_D \v v\v^p d\mu \le A_p \V \mu \V_C  \int_{\bf
T}\v v\v^p dm.\leqno(0.2) $$
We observe in passing that   a simple inner outer
factorisation shows that if (0.1) holds for some $p>0$ then
it also holds for all $p>0$ with the same constant.

It was observed a few years ago   (by J.Bourgain [B], and
also, I believe,  by J.Garcia-Cuerva) that
Carleson's result extends to the Banach space valued case.
More precisely,  there is a numerical constant $K$ such
that, for any Banach space $B$, we have 
  $$\forall p>0 \ \forall f\in H^p(B)\qquad  \int \V f\V^p
d\mu \le K\V \mu \V_C  \V f\V_{H^p(B)}^p .\leqno(0.3)$$ 
Since any
separable Banach  space is isometric to a subspace of
$\l_\i$, this reduces to the following fact. For any
sequence $\{f_n, n\ge 1\}$ in $H^p$, we have $$\int \sup_n
\v f_n\v^p d\mu \le K\V \mu \V_C \int
 \sup_n \v f_n\v^p dm .\leqno (0.4)$$
This can also be deduced from the
scalar case using a simple factorisation argument. Indeed,
let $F$ be the outer function such that $\v F\v = \sup_n \v
f_n\v$ on the circle. Note that by the maximum principle
we have  $\v F\v \ge \sup_n \v f_n\v$ inside $D$, hence
(0.1)   implies
$$\int \sup_n \v f_n\v^p d\mu\le \int   \v F\v^p d\mu\le
C'\int   \v F\v^p dm=\int
 \sup_n \v f_n\v^p dm .$$
This establishes (0.4) (and hence also (0.3)).

We wish to make a connection between Carleson measures
and the following result due to Mireille L\'evy [L]:
\proclaim Theorem 0.1. Let
$S$ be any subspace of $L_1$ and let $u:S\ra L_1(\mu)$ be
an operator. Let $C$ be a fixed constant. Then the
following are equivalent
\item{\rm (i)} For any sequence $\{f_n, n\ge 1\}$ in $S$,
we have $$\int \sup_n\v u(f_n)\v d\mu \le C\int \sup_n\v
f_n\v dm.$$
\item{\rm (ii)} The operator $u$ admits an extension
$\w{u}:L_1\ra L_1(\mu)$ such that $\V \w{u}\V \le C$.

{\bf Proof:\ } This theorem is a consequence of the Hahn
Banach theorem in the same style as in the proof of
theorem 1 below. We merely sketch the proof of
(i)$\Rightarrow$(ii). Assume (i). Let $V\subset
L_\infty(\mu)$ be the linear span of the simple
functions (i.e. a function in $V$ is a linear combination
of disjointly supported indicators). Consider the space 
$S\otimes V$ equipped with the norm induced by the space
$L_1(m;L_\infty(\mu))$. Let $w = \sum\limits^n_1 \varphi_i
\otimes f_i$ with
  $\varphi_i \in V\ \ 
f_i \in S \ $. We will write   $$\a u,w \b =\sum \a
\varphi_i,uf_i\b .$$ Then (i) equivalently means that for
all such $w$ $$|\a u,w \b |\leq C\|
w\|_{L_1(m;L_\infty(\mu))}.$$ By the Hahn Banach theorem,
the linear form $w\ra \a u,w \b$
admits an extension of norm $\leq C$ on the whole of
${L_1(m;L_\infty(\mu))}$. This yields an extension of $u$
from $L_1$ to $L_\infty(\mu)^*=L_1(\mu)^{**}$, with norm $\leq
C$. Finally composing with the classical norm one
projection from $L_1(\mu)^{**}$ to $L_1(\mu) $, we obtain
(ii).

 In
particular, we obtain as a consequence  the following
(known)  fact which we wish to emphasize for later use.
\proclaim Proposition 0.2. Let $\mu$ be a Carleson measure
on $D$, then there is a bounded operator $T:L_1\ra
L_1(\mu)$ such that $T(e^{int})=z^n$ for all $n\ge 0$, or
equivalently such that  $T$ induces the identity on $H^1$.

\n {\bf Proof:} We simply apply L\'evy's theorem to
$H^1$ viewed as a subspace of $L_1$, and to the operator
$u:H^1\ra L_1(\mu)$ defined by $u(f)=f$. By (0.1) we have 
$\V u\V \le K\V \mu \V_C$, but moreover by (0.4) and 
L\'evy's theorem there is an operator $T:L_1\ra
L_1(\mu)$ extending $u$ and with $\V T\V \le  K\V \mu
\V_C$.
This proves the proposition.

Allthough we have not seen this proposition stated
explicitly,  it is undoubtedly known to
specialists (see the remarks below on the operator
$T^*$). Of course, for $p>1$ there is no problem, since in
that case the inequality (0.2) shows that the operator of
harmonic extension (given by the Poisson integral) is
bounded from $L_p$ into $L_p(\mu)$ and of course it
induces the identity on $H^p$. However this same operator
is well known to be unbounded if $p=1$. The adjoint of the
operator $T$ appearing in Proposition 0.2 solves the
$\d$-bar equation in the sense that for  any $\varphi$ in
$L_\i (\mu)$ the function $G=T^*(\varphi)$ satisfies $\V
G\V_{L_\i (m)}\le \V T\V \V \varphi\V_\i$ together with
$$\forall f \in H^1 \quad \int f G  dm =\int
f \varphi d\mu ,$$ and by well known ideas of H\"ormander [H]
this means  equivalently that  $Gdm$ is the boundary value
(in the sense of [H]) of a distribution $g$ on $\bar{D}$
such that $\d g =\varphi .\mu$. In conclusion, we have 
$$\d g
=\varphi .\mu \quad \hbox{and}\quad  \V G \V_{L_\i (m)}\le K\V
\mu \V_C   \V\varphi\V_\i .$$
This is precisely the basic $L_\i$-estimate for the
$\d$-equation proved by Carleson to solve the corona
problem, (cf. [G], theorem 8.1.1, p.320).
More recently, P.Jones [J] proved a refinement of this
result by producing an explicit kernel which plays the
role of the operator $T^*$ in the above. He proved that
one can produce a solution $g$ of the equation $\d g =\varphi
.\mu$ which depends linearly
on $\varphi$ with a boundary value $G$  satisfying 
simultaneously $$\V G \V_{L_\i (m)}\le K\V
\mu \V_C   \V\varphi\V_\i \quad \hbox{and}\quad 
\V G \V_{L_1 (m)}\le K\int \v \varphi\v d\mu
  ,$$ where $K$ is a numerical constant. (Jones [J]
mentions that A.Uchiyama found a different proof of this.
A similar proof, using weights, was later found by
S.Semmes.) Taking into account the previous remarks, our
theorem 1 below gives at the same time a different proof
and a generalization of this theorem of Jones. 

Our
previous paper [P] contains simple direct proofs of
several consequences of Jones' result for interpolation
spaces between $H^p$-spaces. We will use similar ideas in
this paper.

Let us recall here the definition of the
$K_t$   functional which is fundamental for the real
interpolation method.
Let $A_0,A_1$ be a compatible couple of Banach (or
quasi-Banach) spaces. For all $x\in A_0+A_1$ and
for all $t>0$, we let $$K_t(x,A_0,A_1)= \inf
\big({\|x_0\|_{A_0}+t\|x_1\|_{A_1}\ | \ x=x_0+x_1,x_0\in
A_0,x_1\in A_1}).$$
Let $S_0\subset A_0, S_1\subset A_1$ be closed subspaces.
As in [P], we will say that the couple $(S_0, S_1)$ is
$K$-closed (relative to $(A_0, A_1)$) if there is a
constant $C$ such that
$$ \forall t>0 \quad \forall x\in S_0 + S_1 \quad
K_t(x, S_0, S_1) \le C K_t(x, A_0, A_1).$$
\vskip6pt
\n{\bf Main results}
\medskip

\proclaim Theorem 1. Let $(\Omega,\cal{A},\mu)$ be an
arbitrary measure space. Let $u\colon \ H^\infty \to
L_\infty(\mu)$ be a bounded operator with norm $\|u\| =
C_\infty$. Assume that $u$ is also bounded as an operator
from $H^1$ into $L_1(\mu)$, moreover assume that there is
a constant $C_1$ such that for all finite sequences
$x_1,\ldots, x_n$ in $H^1$ we have $$\int \sup_i|u(x_i)|
d\mu \le C_1 \int \sup|x_i|dm.$$ Then there is an operator
$\tilde u\colon \ L^\infty \to L_\infty(\mu)$ which is
also bounded from $L^1$ into $L_1(\mu)$ such that
$$\eqalign{&\|\tilde u\colon \ L_\infty \to
L_\infty(\mu)\|\le CC_\infty\cr &\|\tilde u\colon \ L_1\to
L_1(\mu)\|\le CC_1}$$ where $C$ is a numerical constant.
\medskip

\n {\bf Proof:} \ Let $w$ be arbitrary in $L_\infty(\mu)\otimes H^\infty$.
We introduce on $L_\infty(\mu) \otimes L^\infty(m)$ the following two norms
$\forall w \in L_\infty(\mu)\otimes L^\infty(m)$

$$\eqalign{&\|w\|_0 = \int\|w(\omega,
\cdot)\|_{L^\infty(dm)}d\mu(\omega)\cr
&\|w\|_1 = \int\|w(\cdot, t)\|_{L^\infty(d\mu)}dm(t).}$$

\n Let $A_0$ and $A_1$ be the completions of $L_\infty(\mu)\otimes
L^\infty(m)$ for these two norms. (Note that $A_0$ and
$A_1$ are nothing but respectively $L_{1}(d\mu
;L_{\i}(dm))$ and $L_{1}(dm ;L_{\i}(d\mu ))$. ) 
 Let $S_0$ and
$S_1$ be the closures of $L_\infty(\mu)\otimes H^\infty$
in $A_0$ and $A_1$ respectively.

The completion of the proof is an easy aplication (via
the Hahn Banach theorem) of the following result which is
proved further below:\medskip

\proclaim Lemma  2. $(S_0, S_1)$ is $K$-closed.

Indeed, assuming the lemma
proved for the moment, fix $t>0$, and consider $w$ in
$L_\infty(\mu)\otimes H^\infty$, we have (for some
numerical constant $C$)

$$\forall t > 0 \quad K_t(w, S_0, S_1) \le C K_t(w,
A_0, A_1).$$
Recall that we denote by $V\subset L_\infty(\mu)$ the
dense subspace of functions taking only finitely many
values.  Let $w = \sum\limits^n_1 \varphi_i \otimes f_i$
with
  $\varphi_i \in V\ \ 
f_i \in H^\infty \ $. We will write, for every operator
$u:H^1\ra L_1(\mu)$, $$\a u,w \b =\sum \a \varphi_i,uf_i\b
.$$

Clearly

$$\leqalignno{&\left\|\sum \varphi_i\otimes
u(f_i)\right\|_{L^1_\mu(L^\infty_m)} \le C_\infty
\|w\|_0&(1)\cr \noalign{\hbox{and}}
&\left\| \sum \varphi_i \otimes
u(f_i)\right\|_{L^1_m(L^\infty_\mu)} \le C_1
\|w\|_1&(2)}$$

\n Moreover, by completion, we can extend (1) (resp.
(2)) to the case when $w$ is in $S_0$ (resp. $S_1$).
Hence, if $w = w_0 + w_1$ with $w_0 \in S_0, w_1 \in S_1$
we have by (1) and (2)

$$\eqalign{|\a u,w \b |= \left|\sum\langle \varphi_i,
u(f_i)\rangle\right| &\le C_\infty\|w_0\|_0 +
C_1\|w_1\|_1\cr &\le C_\infty K_s(w, S_0, S_1)\cr
&\le CC_\infty K_s(w, A_0, A_1)}$$

\n where $s = C_1(C_\infty)^{-1}$. By Hahn-Banach, there
is a linear form $\xi$ on $A_0 + A_1$  such that

$$\leqalignno{\xi(w) &= \a u,w \b \qquad \forall w \in S_0 +
S_1&(3)\cr
\noalign{\hbox{and}}
|\xi(w)| &\le CC_\infty K_s(w, A_0, A_1)\qquad \forall w \in
A_0 + A_1.}$$

\n Clearly this implies

$$\leqalignno{\forall w \in A_0\qquad |\xi(w)| &\le
CC_\infty\|w\|_0 &(4)\cr
\forall w \in A_1\qquad |\xi(w)|&\le CC_\infty s\|w\|_1\cr
 &\le CC_1 \|w\|_1.&(5)}$$

\n Now (4) implies $\forall\varphi \in L_\infty(\mu)\quad \forall f \in
L_\infty(dm)$

$$|\langle\xi, \varphi\otimes f\rangle| \le
CC_\infty\|\varphi\|_1 \|f\|_\infty. \leqno (6)$$

\n Define $\tilde u\colon L_\infty \to L_1(\mu)^* = L_\infty(\mu)$ as
$\langle \tilde u(f), \varphi\rangle = \langle \xi, \varphi
\otimes f\rangle$ then (6) implies $\|\tilde
u(f)\|_{L_\infty(\mu)} \le CC_\infty\|f\|_\infty$, while
(5) implies

$$|\langle \xi, \varphi \otimes f\rangle|\le
CC_1\|\varphi\|_\infty \|f\|_1,$$

\n hence $\|\tilde u(f)\|_1\le CC_1\|f\|_1$. Finally (3) implies that
$\forall f\in H^\infty \quad \forall \varphi \in L^\infty(\mu)$

$$\langle \tilde u(f), \varphi\rangle = \langle\varphi, u(f)\rangle$$

\n so that $\tilde u|_{H^\infty}  = u$.\hfill
q.e.d.\bigskip \bigskip

\n {\bf Proof of Lemma  2:} \ We start by reducing this
lemma to the case when $\Omega$ is a finite set or
equivalently, in case the $\sigma$-algebra ${\cal A}$ is
generated by finitely many atoms, with a fixed constant
independent of the number of atoms. Indeed, let $V$ be
  the union of all spaces $L_\infty(\Omega,{\cal
B}, \mu)$ over all the subalgebras ${\cal B}\subset {\cal
A}$ which are generated by finitely many atoms.
Assume the lemma known in that case with a fixed constant
$C$  independent of the number of atoms. It follows
that for any $w$ in $H^\infty \otimes V$ we have

$$\forall t>0\qquad K_t(w, S_0, S_1) \le CK_t(w, A_0,
A_1).$$

\n Since $H^\infty\otimes V$ is dense in $S_0 + S_1$, this
is enough to imply Lemma 2.

Now, if $(\Omega, {\cal B}, \mu)$ is finitely atomic as
 above we argue exactly as
in section 1 in [P] using the simple
(so-called) ``square/dual/square'' argument, as formalized
in Lemma 3.2 in [P]. We want to treat by the same argument
the pair

$$\eqalign{&H^1(L_\infty(\mu)) \subset L^1(L_\infty(\mu))\cr
&L_1(\mu; H^\infty) \subset L_1(\mu; L^\infty).}$$

\n Taking square roots, the problem reduces to prove the following couple
if $K$-closed:

$$\eqalign{&H^2(L_\infty(\mu)) \subset L^2(L_\infty(\mu))\cr
&L_2(\mu; H^\infty)\subset L_2(\mu; L_\infty)}$$

\n provided we can check that

$$H^2(L_\infty(\mu))\cdot L_2(\mu; H^\infty)\subset (H^1(L_\infty(\mu)),
L_1(\mu; H^\infty))_{{1\over 2}\infty}\leqno (7)$$

\n We will check this auxiliary fact below. By duality and by
Proposition 0.1 in [P] , we can reduce to checking the
$K$-closedness for the couple

$$\eqalign{&H^2(L_1(\mu)) \subset L^2(L_1(\mu))\cr
&L_2(\mu; H^1) \subset L_2(\mu; L_1).}$$

\n Taking square roots one more time this reduces to prove that the
following couple is $K$-closed

$$\left\{\eqalign{&H^4(L_2(\mu)) \subset L^4(L_2(\mu))\cr
&L_4(\mu; H^2) \subset L_4(\mu; L_2)}\right.$$

\n provided we have

$$H^4(L_2(\mu))\cdot L_4(\mu; H^2)\subset (H^2(L_1(\mu)), L_2(\mu;
H^1))_{{1\over 2}\infty}. \leqno (8)$$

\n But this last couple is trivially $K$-closed (with a fixed constant
independent of $(\Omega, {\cal B}, \mu)$) because, by Marcel Riesz' theorem,
there is a simultaneously bounded projection

$$\eqalign{&L_4(L_2(\mu))\to H^4(L_2(\mu))\cr
&L_4(\mu; L_2) \to L_4(\mu; H^2).}$$

\n It remains to check the inclusions (7) and (8). We
first check (7). By Jones' theorem (see the
beginning of section 3 and Remark 1.12 in [P])

$$H^2(L_\infty(\mu)) = (H^1(L_\infty(\mu)),
H^\infty(L_\infty(\mu)))_{{1\over 2}2}\leqno (9)$$

\n also by an entirely classical result (cf.[BL] p.109)

$$L_2(\mu; H^\infty) = (L_\infty(\mu; H^\infty), L_1(\mu;
H^\infty))_{{1\over 2}2}.\leqno (10)$$

\n By the bilinear interpolation theorem (cf. [BL]
p.76) the two obvious
inclusions

$$\eqalign{&H^1(L_\infty(\mu))\cdot L_\infty(\mu;
H^\infty) \subset H^1(L_\infty(\mu))\cr
&H^\infty(L_\infty(\mu))\cdot L_1(\mu; H^\infty)\subset L_1(\mu;
H^\infty),}$$

\n (note that $H^\infty(L_\infty(\mu)) = L_\infty(\mu;
H^\infty)$), imply that$$(H^1(L_\infty(\mu)),
H^\infty(L_\infty(\mu)))_{{1\over 2}2}\cdot (L_\infty(\mu; H^\infty), L_1(\mu;
H^\infty))_{{1\over 2}2} \subset (H^1(L_\infty(\mu)),
L_1(\mu; H^\infty))_{{1\over 2}\infty}.$$
Therefore, by (9) and (10), this proves (7).
We now check (8). We will first prove an analogous result
but with the inverses of all indices translated by $1/r$.
More precisely, let $2<r<\i$, let $p,r'$ be defined by
the relations  
$1/2=1/r+1/p$ and $1=1/r+1/{r'}$.
We will first check

 $$H^{2p}(L_{2{r'}}(\mu))\cdot
L_{{2p}}(\mu; H^{2{r'}})\subset (H^p(L_{r'}(\mu)), L_p(\mu;
H^{r'}))_{{1\over 2}\infty}. \leqno (11)$$
Indeed, we have 
$$H^{2p}(L_{2{r'}}(\mu))\cdot
L_{{2p}}(\mu; H^{2{r'}})\subset
L^{2p}(L_{2{r'}}(\mu))\cdot L_{{2p}}(\mu;
L^{2{r'}})  \subset (L^p(L_{r'}(\mu)) , L_p(\mu;
L^{r'}))_{{1\over 2} }. \leqno (12)$$
The last inclusion follows from a classical result on the
complex interpolation of Banach lattices, (cf. [C] p.125).
But now, since all indices appearing are   between $1$
and infinity, the orthogonal projection from $L_2$ onto
$H^2$ defines an operator bounded simultaneously from 
$L^p(L_{r'}(\mu))$ into $H^p(L_{r'}(\mu))$ and from
$L_p(\mu; L^{r'})$ into $L_p(\mu; H^{r'})$, hence also
bounded from   $(L^p(L_{r'}(\mu)) , L_p(\mu;
L^{r'}))_{{1\over 2} }$ into $(H^p(L_{r'}(\mu)) , L_p(\mu;
H^{r'}))_{{1\over 2} }$. Since the latter space is
included into $(H^p(L_{r'}(\mu)) , L_p(\mu;
H^{r'}))_{{1\over 2},\i }$, (cf.[BL] p.102) we obtain the
announced result (11).

 Then, we
use the easy fact that any element $g$ in the unit
ball of $H^4(L_2(\mu))$(resp. $h$ in the unit
ball of $ L_4(\mu; H^2)$) can be
written as $g=Gg_1$ (resp. $h=Hh_1$) with $G$ and $H$ in
the unit ball of $H^{2r}(L_{2r}(\mu))=L_{2r}(\mu;H^{2r})$
and with $g_1$ (resp. $h_1$) in the unit ball of
$H^{2p}(L_{2{r'}}(\mu))$ (resp. $ L_{{2p}}(\mu;
H^{2{r'}})$). Then, by (11), there is a constant $C$
such that
$$\V g_1h_1\V_{(H^p(L_{r'}(\mu)), L_p(\mu;
H^{r'}))_{{1\over 2}\infty}} \le C.$$ 
Now, the product $M=GH$ is in the unit ball of
$H^{{r}}(L_{{r}}(\mu))=L_{{r}}(\mu;H^{{r}})$, therefore
the operator of multiplication by $M$ is of norm $1$ both
from   $H^p(L_{r'}(\mu))$ into $H^2(L_{1}(\mu))$ and
from $L_p(\mu;
H^{r'})$ into $L_2(\mu;
H^{1})$. By interpolation, multiplication by $M$ also has
norm $1$ from ${(H^p(L_{r'}(\mu)), L_p(\mu;
H^{r'}))_{{1\over 2}\infty}}$ into ${(H^2(L_{1}(\mu)),
L_2(\mu; H^{1}))_{{1\over 2}\infty}}$. Hence, we conclude
that $gh=Mg_1h_1$ has norm at most $C$ in the space ${(H^2(L_{1}(\mu)),
L_2(\mu; H^{1}))_{{1\over 2}\infty}}$. This concludes the
proof of (8).

\centerline {\bf References}\vskip6pt

\item {[BL]} J.Bergh and J.L\"ofstr\"om, Interpolation
spaces, An introduction, Springer Verlag 1976.

\item {[BS]} C.Bennett and R.Sharpley, Interpolation of
operators.Academic Press,1988.

\item {[B]} J.Bourgain. On the similarity problem for
polynomially bounded operators on Hilbert space, Israel
J. Math. 54 (1986) 227-241.
 
\item {[C]} A.Calder\'on, Intermediate spaces and
interpolation, Studia Math. 24 (1964) 113-190.

\item {[G]} J.Garnett, Bounded Analytic Functions.
Academic Press 1981.

\item {[GR]} J.Garcia-Cuerva and J.L.Rubio de Francia.
Weighted norm inequalities and related topics. North
Holland, 1985.

\item {[H]} L.H\"ormander, Generators for some rings of
analytic functions, Bull.Amer.Math.Soc. 73 (1967) 943-949.

\item {[J]} P.Jones, $L^\infty$ estimates for the
$\bar{\partial}$-problem in a half plane. Acta Math. 150
(1983) \nobreak{137-152}.

\item {[L]} M.L\'evy. Prolongement d'un op\'erateur d'un
sous-espace de $L^1(\mu)$ dans $L^1(\nu)$. S\'emi-naire
d'Analyse Fonctionnelle 1979-1980. Expos\'e 5. Ecole
Polytechnique.Palaiseau.

 \item {[P]} G.Pisier, Interpolation between $H^p$ spaces
and non-commutative generalizations I. Pacific J. Math.
(1992) To appear.

\vskip12pt

Texas A. and M. University

College Station, TX 77843, U. S. A.

and

Universit\'e Paris 6
Equipe d'Analyse, Bo\^{\A}te 186,
4 Place Jussieu, 75230 

Paris Cedex 05, France

\end

\item {[X1]} Q.Xu, Applications du th\'eor\`eme de
factorisation pour des fonctions \`a valeurs op\'e
-rateurs.  Studia Math. 95 (1989) 273-292.

\item {[X2]} Q.Xu, Real interpolation of some Banach
lattices valued Hardy spaces. Bull. Sci. Mat.
  116 (1992) To appear.

\item {[X3]} Q.Xu, Elementary proofs of two theorems of
P.W.Jones on interpolation of Hardy spaces. Preprint,  
Pub.\quad Irma,
Lille, 1989.

\item {[HP]} U. Haagerup and G. Pisier, Factorization of
analytic functions with values in non-commutative
$L_1$-spaces and applications. Canadian J. Math. 41
(1989) 882-906.

\item {[J2]} P.Jones, Interpolation between Hardy
spaces,  in: Conference on Harmonic Analysis in
honor of Antoni Zygmund, (edited by W.Beckner,
A.Calder\'on,R.Fefferman and P.Jones) 
 Wadsworth Inc.,1983, vol.2, p.437-451.

\item {[J3]} P. Jones. On interpolation
between $H^1$ and
$H^\infty$ (in Interpolation Spaces and Allied Topics in Analysis) Springer
L.N. no 1070 (1984) 143-151\bigskip

\item {[JJ]} S.Janson and P.Jones, Interpolation between
$H^p$ spaces: The complex method.
Journal Funct. Anal. 48 (1982) 58-80.
 H.Helson, Lectures on invariant
subspaces.Academic Press, New-York 1964.
\item {[GK]} I.C.Gohberg and M.G.Krein, Introduction
to the theory of linear nonselfadjoint
operators, Transl. Math. Monogrphs, Amer. Math. Soc.
Providence, RI, 1969.  

\item {[W]} T.Wolff. A note on interpolation spaces.
Springer Lecture Notes in Math. 908 (1982) 199-204.